\documentclass[12pt]{amsart}
\usepackage{amssymb}
\newfont{\ams}{msbm10 at 12pt}
\newfont{\amsi}{msbm8}
\newcommand{\fg}{{\mathfrak g} }

\newcommand{\fL}{{\mathfrak L} }
\newcommand{\cC}{{\mathcal C} }
\newcommand{\cA}{{\mathcal A} }
\newcommand{\oC}{\overline{C} }
\newcommand{\ow}{\overline{w} }
\newcommand{\cQ}{{\mathcal Q} }
\newcommand{\cK}{{\mathcal K} }
\newcommand{\cO}{{\mathcal O} }
\newcommand{\cL}{{\mathcal L} }
\newcommand{\cN}{{\mathcal N} }
\newcommand{\cH}{{\mathcal H} }
\newcommand{\BZ}{{\mathbb Z} }

\newcommand{\BC}{{\mathbb C} }
\newcommand{\BN}{{\mathbb N} }
\begin{document}
\title{Tensor ideals in the category of tilting modules.}
\authors{V.Ostrik\thanks{ This material is based upon work supported
by the U.S. Civilian
Research and Development Foundation under Award No. RM1-265.}
\address{Independent Moscow University, 11 Bolshoj Vlasjevskij per.,
 Moscow
121002 Russia}
\date{November 1996}
\email{ostrik@nw.math.msu.su}
\maketitle

\section{Introduction.}
Let $\fg$ be a complex finite dimensional simple Lie algebra with the
 root
datum $(Y,X,\ldots ),$ see \cite{Lu2}.
 Let $W_f$
denote the Weyl group, $R$ denote the root system, $R_+$ denote the
 set
 of positive roots.
Let $X_+$ denote the set of dominant integral weights.
Let $h$ denote the Coxeter number of $\fg.$

Let us fix $l\in \BN,\, l>h.$ We assume that $l$ is odd (and not
 divisible by
 3,
if $\fg$ is of type $G_2$).
 Let $W$ denote the corresponding affine Weyl group.

Let $\rho \in X$ denote the halfsum of positive roots.
 We will denote by dot
(for example $w\cdot \lambda$) the action of $W$ (and $W_f\subset W$)
 centered in $(-\rho).$

 Let $q$
be a primitive $l-$th root of unity and let $U_q$ be the quantum
group with divided powers as defined in \cite{Lu2}.
 Let $\cC$ denote the category of finite dimensional $U_q-$modules of
type {\bf 1} (see e.g. \cite{APW}).

In \cite{An1} H.Andersen has studied a tensor subcategory
 $\cQ \subset \cC$
formed by {\em tilting} modules. He has introduced a tensor ideal
 $\cK \subset \cQ$ formed by negligible tilting modules. The quotient
tensor category $\cQ /\cK $ is semisimple. For certain values of $l$
 it is
tensor-equivalent to a category of integrable modules over affine Lie
algebra $\hat \fg$ equipped with a {\em fusion} tensor structure
(see e.g. \cite{F}).

Let us recall the definition of $\cK.$ Indecomposable tilting modules
 are
numbered by their highest weights $\lambda \in X_+;$ we will denote
 them
by $Q(\lambda).$ The set of dominant weights $X_+$ is covered by the
 closed
{\em alcoves} numbered by $W^f\subset W$ --- the set of shortest
 elements
in the right cosets $W/W_f.$ For $w\in W^f$ the corresponding closed
 alcove
will be denoted by $\oC_w.$ For example, the alcove $\oC_e=\oC$
 containing
the zero weight is given by
$$
\oC=\{ \lambda\in X|0\le \langle \lambda +\rho
 ,\alpha^{\vee}\rangle \le
 l \, \mbox{for all}\, \alpha \in R_+\}.
$$
Now $\cK$ is formed by the direct sums of tiltings $Q(\lambda),$
where $\lambda$ is dominant and
$\lambda \in \bigcup_{w\ne e}\oC_w.$

In this note we propose the following generalization of H.Andersen's
 result.
We recall that G.Lusztig and N.Xi have introduced a partition of
 $W^f$ into
{\em canonical right cells} along with the {\em right order} $\le_R$
 on the
set of cells, see \cite{Lu1} and \cite{LX}. In particular, $\{ e\}
 \subset
W^f$
forms a single right cell, maximal with respect to $\le_R.$ Thus
$W^f-\{ e\} =\coprod_{A<_R\{ e\} }A$ --- the union of right cells.

{\bf Main Theorem.} {\em Let $A\subset W^f$ be a right cell. The
full subcategory
$\cQ_{\le A}$ formed by the direct sums of tiltings $Q(\lambda),\;
\lambda \in \bigcup_{w\in B\le_RA}\oC_w,$ is a tensor ideal
 in $\cQ.$}

There is a well-known correspondence between the right cells in $W$
 and the
right ideals in the affine Hecke algebra $\cH$ (see \cite{KL}). Our
 result
is completely parallel to this correspondence, and even the proof is.
In fact, the proof is an application of a deep result by W.Soergel
 who
has connected the characters of $Q(\lambda)$ with
 Kazhdan-Lusztig-type
combinatorics of $\cH.$

In general, the right cells in $W^f$ are infinite, but some are
 finite, e.g.
$\{ e\} \subset W^f.$ The first nontrivial example is a "subregular"
 cell
$D_1$ for $\fg$ of type $G_2$ (see the pictures and notations
 in \cite{Lu1}) consisting
of 8 alcoves. Then the subcategory $\cQ_{<D_1}$ formed by the direct
 sums
of $Q(\lambda)$ such that $\lambda \in \bigcup_{w\in B<_RD_1}\oC_w$
 is a
tensor ideal, and we can consider the quotient subcategory
 $\cQ /\cQ_{<D_1}$
with finitely many isomorphism classes of indecomposable objects.
 This
subcategory is nonsemisimple, as opposed to Andersen's fusion
 category
$\cQ /\cK.$ For example, when $l=7,\; \cQ/\cK$ is equivalent to
 $\BC -$vector
spaces, while $\cQ/\cQ_{<D_1}$ has 24 isomorphism classes of
 indecomposable
objects. Its Grothendieck ring is a 24-dimensional
 algebra with nontrivial nilpotent radical, as
opposed to the classical fusion rings which are always
 semisimple. To our
knowledge, this is a first example of a nonsemisimple
 tensor category
without fiber functor with finitely many indecomposable
 objects.

As we already mentioned, for certain values of $l,$ $\cQ/\cK$
 is  tensor
equivalent to a category of integrable $\hat \fg -$modules of
 positive
central charge. It is a subcategory of a larger category $\cO$ of all
$\fg -$integrable $\hat \fg -$modules of positive central charge,
 but the
Kazhdan-Lusztig construction of fusion tensor structure in this
 larger
 category encounters serious problems (see \cite{KL2}). Still
 we believe
that the quotient categories like $\cQ/\cQ_{<D_1}$ are closely
 related to
the would-be fusion structure on $\cO.$

The idea of this note is essentially due to J.Humphreys : it was
 he who
suggested the important role played by the right cells in the
 study of
tilting modules \cite{h}. I learnt of his ideas from M.Finkelberg.
I am grateful to  Catharina Stroppel for her beautiful patterns of
tilting characters for $G_2$ which provided a further insight
 into the
connection between right cells and tilting modules.
Thanks are also due to
D.Timashov who acquainted me with LIE package; it was very useful
 for me
at the first stage of my work.
 I am indebted to H.H.Andersen and J.Humphreys
 for the valuable suggestions
which improved the exposition.
Finally, I would like to thank the referee for extremely
useful comments which simplified the original proof drastically.

 \section{Preliminaries.}

\subsection{}
For any $\lambda \in \overline{C}$ let $\cC (\lambda)$ denote a full
subcategory
of $\cC$ consisting of modules whose composition factors have highest
weights in $W\cdot \lambda.$ The category $\cC$ is a direct sum of
the subcategories $\cC (\lambda)$ (linkage principle; see
 e.g. \cite{APW},
\S 8)
$$
\cC=\bigoplus_{\lambda \in \overline{C}}\cC (\lambda).
$$

For any $\lambda \in X_+$ one defines Weyl module $V(\lambda)$ and
 module
$H^0(\lambda)$ (see \cite{An1} \S 1). Then the irreducible module
 $L(\lambda)$
is the socle of $H^0(\lambda)$ as well as the head of $V(\lambda).$

\subsection{}
\label{char}
Let $\BZ [X]$ be the group algebra of abelian group $X.$
 It is generated by
elements $e^{\lambda},\, \lambda \in X,$ with relations
 $e^{\lambda_1}\cdot e^{\lambda_2}=e^{\lambda_1+\lambda_2} \;
 \forall \lambda_1,\lambda_2 \in X.$
 There is a natural
action of $W_f$ on $\BZ [X]$ given by the formula
 $we^{\lambda}=e^{w\lambda}.$
 Let $\cA :=\BZ [X]^{W_f}$ be the invariants of this action.
 It is a subalgebra
of $\BZ [X].$

Let $ch:K(\cC) \to \BZ [X]$ be the map associating to a module
 $M\in \cC$ its character $ch(M).$ It is known
that its image is $\cA.$ Moreover the elements $ch([V(\lambda)])$
 where
 $\lambda$ runs through $X_+$ form a basis of $\cA.$ It is known
that $ch([V(\lambda)])=ch([H^0(\lambda)])$
is given by the Weyl character formula (see e.g. \cite{APW} \S 8):
$$
ch([V(\lambda)])=\frac{\sum_{w\in W_f}(-1)^{l(w)}e^{w\cdot
\lambda}}{\sum_{w\in W_f}(-1)^{l(w)}e^{w\cdot 0}}.
$$

Now for any $\lambda \in X$ let
$$
ch(\lambda)=\frac{\sum_{w\in W_f}(-1)^{l(w)}e^{w\cdot
\lambda}}{\sum_{w\in W_f}(-1)^{l(w)}e^{w\cdot 0}}.
$$

{\bf Lemma.} (i){\em If stabilizer (in $W_f$) of $\lambda$
with respect to dot action is
nontrivial then $ch(\lambda)=0.$  }

(ii){\em Suppose the stabilizer of $\lambda$ is trivial and let
$w\in W_f$ be  such
that $w\cdot \lambda \in X_+.$ Then $ch(\lambda)=(-1)^{l(w)}
ch(w\cdot \lambda).$}

{\bf Proof.} Clear. $\square$

\subsection{}
\label{calc}
Let $W\to W_f,\, w\mapsto \ow$ be the standard homomorphism with the
 kernel
consisting of translations.

{\bf Lemma.} {\em For any $\lambda,\mu \in X$ and $w\in W$ we have
$w(\lambda+\mu)=w\lambda+\ow \mu$ and
 $w\cdot (\lambda+\mu)=w\cdot \lambda+\ow \mu.$}

{\bf Proof.} The first identity is obviously true for
 $w\in W_f\subset
 W$ and
for translations. Since $W$ is a semidirect product of $W_f$
 and the subgroup of
translations we get our result. The second identity is  a simple
 consequence
of the first one. $\square$

\subsection{}
\label{calc1}
{\bf Lemma.} (see e.g. \cite{D} 2.2.3)
{\em Let $P$ be a multiset (set with multiplicities)
 of weights
invariant under $W_f$ action. Then for any $\lambda \in X$ we have
$$
(\sum_{\omega \in P}e^{\omega})ch(\lambda)=\sum_{\omega \in P}
ch(\lambda+\omega)
$$
}

{\bf Proof.} Straightforward computation. $\square$

\subsection{}
A filtration of $U_q-$module is called {\em Weyl filtration}
 (respectively
{\em good filtration}) if all the associated factors are Weyl
 modules
(respectively modules $H^0(\lambda)$).

\subsection{}
\label{tilting}
{\bf Definition} (see \cite{An1}, definition 2.4) {\em A tilting
 module is
a module
$M\in \cC$ which has both a Weyl filtration and a good filtration.}

Let $\cQ \subset \cC$ be a full subcategory formed by all tilting
 modules.
 The main properties of this
category are collected in the following (see \cite{An1} \S 2)

{\bf Theorem.} (i) {\em The category $\cQ$ is closed under tensor
multiplication.}

 (ii)  {\em
Any tilting module is a sum of indecomposable tilting modules.}

 (iii) {\em For each $\lambda \in X^+$ there exists an indecomposable
tilting module $Q(\lambda)$ with highest weight $\lambda.$ }

 (iv) {\em The modules $Q(\lambda), \lambda \in X^+,$ form
 a complete
set of nonisomorphic indecomposable tilting modules}

 (v) {\em A tilting module is determined up to isomorphism by its
character.}

Let $\cQ(\lambda)$ be the full subcategory of $\cQ$ consisting of
 modules
contained in $\cC(\lambda).$ Then obviously
$$
\cQ=\bigoplus_{\lambda \in \overline{C}}\cQ (\lambda).
$$

\subsection{}
For any $\lambda,\mu \in \oC$ one defines the translation
functor $T_{\lambda}^{\mu}:\cC(\lambda)\to \cC(\mu)$
(see e.g. \cite{APW} \S 8). The following Proposition
is proved as in \cite{J},II,7.13.

\subsubsection{}
{\bf Proposition.} {\em Suppose $\lambda,\mu \in \oC$
and $w\in W$ is such that $w\cdot \lambda \in X_+$.
Then $T_{\lambda}^{\mu}V(w\cdot \lambda)$ has a filtration
with the associated factors $V(\nu)$ such that $\nu \in X_+$
and $\nu=ww_1\cdot \mu$ with $w_1\in Stab(\lambda)$. Each
one of the above factors occurs exactly once.}

In particular it follows that translation functors preserve
the category $\cQ$.

\subsubsection{}
\label{tr}
{\bf Corollary.} {\em For any $w\in W$ such that
$w\cdot \lambda \in X_+$  the module $T^0_{\lambda}T^{\lambda}_0
V(w\cdot 0)$ has a filtration with associated factors $V(wx\cdot 0)$
with $x\in Stab(\lambda)$.}

{\bf Proof.} Evident. $\square$

\section{Construction of tensor ideals.}

\subsection{}
Recall that $W$ denotes the affine Weyl group and $W_f$ denotes
 the ordinary
Weyl group. Let $W^f$ denote the set of minimal length
 representatives of right
cosets. The multiplication defines a bijection
 $W_f\times W^f\to W.$ Let
$\fL$ be the sign representation of $W_f.$ We will consider it
 as right
 $W_f-$module. Let us define a right $W-$module
 $\cN^1:=\fL \otimes _{\BZ [W_f]} \BZ [W].$ As $\BZ -$module it
 is isomorphic
to a free abelian group with generators numbered by $W^f.$ Let
$N^1_x=1\otimes x$ for any $x\in W^f.$ These elements form a
 $\BZ -$basis
of $\cN^1.$ For any $s\in S$ we have $N^1_xs=N^1_{xs}$ if
 $xs\in W^f$ and
$N^1_xs=-N^1_x$ otherwise.

\subsection{}
 Let $K(\cC)$ denote the
Grothendieck group of the category $\cC$.
For any $\lambda \in \oC$ define the map
$\alpha_{\lambda}:K(\cC)\to \cN^1$ by
$\alpha_{\lambda}([V(\mu)])=1\otimes
( \sum_{x\in W, x\cdot \lambda=\mu}x).$
In particular $\alpha_{\lambda}$ annihilates every object
outside of the block $\cC(\lambda)$ of $\cC$.

\subsubsection{}
Let us identify $K(\cC)$ with the character ring $\cA$.

{\bf Lemma.} {\em For any $w\in W$ we have
$$
\alpha_{\lambda}(ch(w\cdot \lambda))=|Stab(\lambda)|^{-1}1\otimes
 (\sum_{x\in Stab(\lambda)}wx)
$$
}

{\bf Proof.} For $w\in W$ such that $w\cdot \lambda \in X_+$ the
Lemma
is clear from definitions. For other $w$ use \ref{char}.  $\square$

\subsubsection{}
{\bf Lemma.} {\em For any $\lambda \in \oC$
and $V \in \cC (\lambda)$ we have
$$
\alpha_{\lambda}(V)=\alpha_0(T_{\lambda}^0V)
$$ }

{\bf Proof.} Obvious. $\square$

\subsubsection{}
\label{last}
{\bf Lemma.} {\em For any $\lambda \in \oC$
and $V \in \cC (0)$ we have
$$
\alpha_0(T_{\lambda}^0T_0^{\lambda}V)=
\alpha_0(V)\sum_{x\in Stab(\lambda)}x
$$ }

{\bf Proof.} It is enough to verify the Lemma
for $V=V(w\cdot 0)$. Now if $w\cdot \lambda \in X_+$
the result follows from \ref{tr}; if $w\cdot \lambda \not \in X_+$
then RHS and LHS both vanish. $\square$

\subsection{}
{\bf Proposition.} {\em For any $\lambda,\mu \in \oC$ and
$M\in \cC$ there exists $c(M)=c_{\lambda \mu}(M)\in \BZ [W]$
such that for all $V\in C(\lambda)$ we have
$$
\alpha_{\mu}(V\otimes M)=\alpha_{\lambda}(V)c(M).
$$}

{\bf Proof.} (see also \cite{J},II,7.5)
 It is enough to check the claim on the level of characters;
moreover we can suppose that $ch(V)=ch(w\cdot \lambda)$.

 Let $P(M)$ be a multiset of weights of module
$M$. It is invariant under $W_f-$action.
 We have by \ref{calc1} and \ref{calc}
$$
ch(V(w\cdot \lambda)\otimes M)=
\sum_{\omega \in P(M)}ch(w\cdot \lambda+\omega)=
\sum_{\omega \in P(M)}ch(w\cdot (\lambda +\omega))
$$

Now let us define a multiset $W_{\lambda \mu}(M):=\{ x\in W|
\lambda +\omega=x\cdot \mu; \omega \in P(M)\} $ It is easy
to see that $W_{\lambda \mu}(M)$ is invariant under left
multiplication by elements of $Stab(\lambda)$ and right
multiplication by elements of $Stab(\mu)$.
So $W_{\lambda \mu}(M)$
is a union of left and right cosets; let $W_{\lambda \mu}(M)'$ be a
set of representatives of right cosets.
We claim that we can choose $c_{\lambda\mu}(M)
=\sum_{z\in W_{\lambda \mu}(M)'}z$.

Indeed,
let $P_{\lambda\mu}(M):=
\{ \omega \in P(M)|\lambda +\omega \in W\cdot \mu \}$.
For any $\omega \in P_{\lambda\mu}(M)$ let $w(\omega)$ be any
 element of $W$ such that
$w(\omega)^{-1}\cdot (\lambda +\omega)=\mu$. It is evident
that $\{ w(\omega)\}$ is the set of representatives of left cosets in
$W_{\lambda\mu}(M)$. We have
$$
\alpha_{\mu}(V\otimes M)=\alpha_{\mu}(ch(w\cdot \lambda)ch(M))=
\alpha_{\mu}(\sum_{\omega \in P(M)}ch(w\cdot (\lambda+\omega)))=
$$
$$
\alpha_{\mu}(\sum_{\omega \in P_{\lambda\mu}(M)}
ch(w\cdot (\lambda+\omega)))=
\alpha_{\mu}(\sum_{\omega \in P_{\lambda\mu}(M)}
ch(ww(\omega)\cdot \mu))=
$$
$$
\sum_{x\in Stab(\mu)}\sum_{\omega \in P_{\lambda\mu}(M)}
1\otimes ww(\omega)x=
 \sum_{t\in W_{\lambda\mu}(M)}1\otimes wt=
$$
$$
\sum_{y\in Stab(\lambda)}\sum_{z\in W_{\lambda\mu}(M)'}
1\otimes wyz=
\alpha_{\lambda}(ch(w\cdot \lambda))\sum_{z\in W_{\lambda\mu}(M)'}z
$$

The Proposition is proved. $\square$

\subsection{}
{\bf Definition.} {\em A subcategory $\cC'\subset \cC$ is called
a weak tensor ideal if for any $V\in \cC'$ and $M\in \cC$ we have
$V\otimes M\in \cC'$.}

We define weak tensor ideals in any subcategory of $\cC$
closed under tensor multiplication in the same way.

{\bf Corollary.} {\em If $U\subset \cN^1$ is a $\BZ [W]-$submodule,
then $\cC_U:=\{V\in \cC|\alpha_{\lambda}(V)\in U \;
 \forall \lambda \in \oC \}$ is a weak tensor ideal of $\cC$ and
$\cQ_U:=\cQ \cap \cC_U$ is a weak tensor ideal of $\cQ$.}

{\bf Proof.} Clear. $\square$

\section{Realization of $K(\cQ(0))$ as a module over Hecke algebra.}
In this section we follow \cite{So1} .

\subsection{}
 Let
$l:W\to \BN$ be the length function and let $\le $ be the standard
Bruhat order on $W.$ We will write $x<y$ if $x\le y$ and $x\ne y.$
Let $\cL=\BZ[v,v^{-1}]$ denote the ring of Laurent polynomials over
 $\BZ$ in
variable $v.$ Let $\cH$ be the Hecke algebra corresponding to $(W,S)$
$$
\cH=\bigoplus_{x\in W} \cL T_x
$$
 with multiplication given by the rule:
 $T_xT_y=T_{xy}$ if $l(xy)=l(x)+l(y)$
and $T_s^2=v^{-2}T_e+(v^{-2}-1)T_s$ for all $s\in S$
 (see \cite{So1} \S 2).

Let $H_x=v^{l(x)}T_x$ be a new basis of Hecke algebra. There exists
 unique
involutive automorphism of Hecke algebra $d:\cH \to \cH, H\mapsto
\overline{H}$
such that $\overline{v}=v^{-1}$ and
 $\overline{H_x}=(H_{x^{-1}})^{-1}.$
We will call $H\in \cH$ selfdual if $\overline{H}=H.$

The following theorem was proved by Kazhdan and Lusztig in \cite{KL}.

{\bf Theorem.}{\em For any $x\in W$ there exists unique selfdual
 $\underline{H}_x\in \cH$ such that
        $\underline{H}_x\in H_x+\sum_{y<x} v\BZ[v]H_y.$}

The coefficients of $\underline{H}_x$ in the basis $\{H_x\}$ are
 essentially
Kazhdan-Lusztig polynomials.

\subsection{}
Let $\cH _f$ be the Hecke algebra corresponding to the group $W_f.$
 We have
an obvious embedding $\cH _f\subset \cH.$ Let $\cL (-v)$ be a free
 right
 $\cL-$module of rank 1 with the right action of $\cH _f$ given by
 the
 following rule:
for any $s\in S_f$ the element $H_s$ acts as $(-v).$ We define a
 right
$\cH -$module $\cN :=\cL (-v)\otimes_{\cH _f}\cH.$
 For any
 $x\in W^f$ let us define
$\underline{N}_x:=1\otimes \underline{H}_x\in \cN.$
 Let
 $\beta :\cN \to \cN^1$
denote the specialization map: $v \mapsto 1.$
We define
$\underline{N}_x^1:=\beta (\underline{N}_x)\in \cN ^1$.

\subsection{}
\label{ch}
The following statement was conjectured in \cite{So1} (Vermutung 7.2)
     and then proved in \cite{So2}.

{\bf Theorem.} {\em
$\alpha (Q(x\cdot 0))=\underline{N}_x^1$.}

\subsection{}
We will say that an $\BZ [W]-$submodule of $\cN ^1$ is a
{\em KL-submodule}
if it admits a base consisiting of elements $\underline{N}_x^1$
 for some subset of $W^f$.

\section{Right cells in affine Weyl group.}

\subsection{}
In \cite{KL} Kazhdan and Lusztig defined three partitions of any
 Coxeter
group into subsets called right, left and two-sided cells
 respectively. We
refer
the reader to {\em loc. cit.} for the definitions of preorders
 $\le_R,\, \le_L, \le_{LR}$
on Coxeter groups. The right (left, two-sided) cells are the
 classes of
equivalence generated by preorder $\le_R$ (respectively $\le_R$ and
$\le_{LR}$).
Let $w\in W$ and $A$ be a right cell in $W.$ We will write that
 $w\le_RA$
if $w\le_Rw'$ for any $w'\in A$ (and similarly for left and
 two-sided cells).

\subsection{}
\label{corr}
There is a correspondence between cells and ideals in the Hecke
 algebra.
 Namely, for any right (left or two-sided)
cell $A$ the $\cL-$submodule $I_{\le A}$ of $\cH$ generated by
 $\underline{H}_w,\,
 w\le_RA$ (and similarly for left and two-sided cells) is a
right (respectively left and two-sided) ideal of $\cH $
 (see \cite{KL}). Moreover any KL-ideal (i.e. ideal admiting
a base consisting of some elements $\underline{H}_w$) is a sum
of such ideals.

\subsection{}
Let $A$ be a two-sided cell of $W.$ The main result of \cite{LX}
 is the
following

{\bf Theorem.} {\em The intersection $A\cap W^f$ forms a right cell
of $W.$ }

\subsection{}
 {\bf Definition.} {\em A weak tensor ideal $\tau \subset \cQ$
 is called a tensor
ideal if for any $Q_1, Q_2$ such that $Q_1\oplus Q_2 \in \tau$
we have $Q_1, Q_2 \in \tau$.}

For any two-sided cell $A$ of $W$ we define the full subcategory
 $\cQ_{\le A}$
 of $\cQ$ as follows: $\cQ_{\le A}$ is the additive subcategory of
 $\cQ$
and indecomposable objects of $\cQ_{\le A}$ are all the modules
$Q(w\cdot \lambda)$ where $\lambda \in \oC ,\, w\in W^f$ and
 $w\le_RA.$

\subsection{Main Theorem}
{\em For any two-sided cell $A$ of $W$ the subcategory $\cQ_{\le A}$
 is a
tensor ideal.}

{\bf Proof.} For any two-sided cell $A$ we define a
$\BZ [W]-$submodule
$U_{\le A}$ of $\cN^1$ to be $\cL \otimes I_{\le A\cap W^f}$.

 We will show that for any $\lambda \in \oC \; \;$
 $\alpha_{\lambda}
(Q(w\cdot \lambda))\in U_{\le A}$ if and only if
$\alpha_0(Q(w'\cdot 0) \in U_{\le A}$ where $w'$ is the longest
element of coset $wStab(\lambda)$. We have
$$
\alpha_{\lambda}(Q(w\cdot \lambda))=|Stab(\lambda)|^{-1}
\alpha_0(T_{\lambda}^0Q(w\cdot \lambda))
$$
Note that $T_{\lambda}^0Q(w\cdot \lambda)$ contains a direct
summand $Q(w'\cdot 0)$. So we proved that $\alpha_{\lambda}
(Q(w\cdot \lambda))\in U_{\le A}$ implies that
$\alpha_0(Q(w'\cdot 0) \in U_{\le A}$.

Now note that $T_0^{\lambda}Q(w'\cdot 0)$ contains a direct
summand $Q(w'\cdot \lambda)=Q(w\cdot \lambda)$. Further
$\alpha_0(T_{\lambda}^0T_0^{\lambda}Q(w'\cdot 0))=
\alpha_0(Q(w'\cdot 0))\sum_{x\in Stab(\lambda)}x\in U_{\le A}$
by \ref{last} and we proved our claim in another direction.

So the proof of theorem is finished. $\square$

\subsection{}
{\em Remark.} It is easy to see that theorem above
establishes bijection between KL-submodules of $\cN ^1$
and tensor ideals in $\cQ$.
Further note that all KL-submodules of $\cN^1$ are the
sums of submodules $U_{\le A}$. So we describe all
tensor ideals in a category of tilting modules.


\begin{thebibliography}{99}

\bibitem[A]{An1} H.H.Andersen {\em Tensor products of quantized
 tilting
modules.} Communications in Mathematical
 Physics 149 (1992), pp. 149-159.

\bibitem[APW]{APW} H.Andersen, P.Polo, K.Wen {\em Representations of
quantum algebras.} Invent. Math. 104 (1991), pp.1-59.

\bibitem[D]{D} S.Donkin
{\em Rational representations of algebraic groups.} Lecture Notes
in Math. vol 1140, Berlin Heidelberg, New York: Springer 1985.

\bibitem[F]{F} M.Finkelberg
 {\em An equivalence of fusion categories},
Geometrical and Functional Analysis, vol.6, 2 (1996), pp. 249-267.

\bibitem[H]{h} J.E. Humphreys, {\em Comparing modular representations
 of
semisimple groups and their Lie algebras}, in volume dedicated to
 R.E. Block,
Modular Interfaces: Modular Lie Algebras, Quantum Groups, and Lie
Superalgebras, ed. V. Chari and I. Penkov, International Press, 1997.

\bibitem[J]{J} J.C.Jantzen {\em Representations of Algebraic Groups},
 Pure
and Applied Mathematics 131, Academic Press, 1987.

\bibitem[KL1]{KL} D.Kazhdan, G.Lusztig {\em Representations of
 Coxeter groups
and Hecke algebras} Inventiones math. 53 (1979), pp. 165-184.

\bibitem[KL2]{KL2} D.Kazhdan, G.Lusztig {\em Tensor structures
 arising
from affine Lie algebras I-IV}, J.Amer.Math.Soc. 6 (1993),
 pp. 905-1011 and
J.Amer.Math.Soc. 7 (1994), pp. 335-453.

\bibitem[L1]{Lu1} G.Lusztig {\em Cells in affine Weyl groups.},
 Advanced
Studies in Pure Mathematics 6, 1985; Algebraic Groups and Related
Topics pp. 255-287.

\bibitem[L2]{Lu2} G.Lusztig {\em Introduction to quantum groups.},
 Boston,
Birkhauser, 1993.

\bibitem[LX]{LX} G.Lusztig, N.Xi {\em Canonical left cells in
 affine Weyl
groups,} Advances in Mathematics 72 (1988), pp. 284-288.

\bibitem[S1]{So1} W.Soergel {\em Kazhdan-Lusztig-Polynome und eine
Kombinatorik f\"ur Kipp-Moduln.}, Representation Theory, 1 (1997).

\bibitem[S2]{So2} W.Soergel {\em Charakterformeln f\"ur Kipp-Moduln
\"uber Kac-Moody-Algebren.}, Preprint (1996), pp. 1-21.

\end{thebibliography}
\end{document}